\documentclass[a4paper 11pt]{article}
\usepackage{latexsym,amsmath, amsthm, amssymb}
\catcode`@=11 \long\def\@makefntext#1{\noindent #1}
\newskip\tabcentering \tabcentering=1000pt plus 1000pt minus 1000pt
\def\MCH#1#2{\setbox0=\hbox{\raise#1\hbox{#2}}\smash{\box0}}% move char

\def\@evenfoot{}\def\@oddfoot{}

\def\@evenhead{\hbox to\textwidth{\footnotesize\rm\thepage \hfill
{\it }}}

\def\@oddhead{\hbox to \textwidth{\footnotesize{\it
} \hfill\thepage}}

\def\proof{\vspace{2mm}\noindent{\it Proof}\quad}

\newtheorem{thm}{Theorem}[section]
\newtheorem{definition}{Dfinition}[section]

\newtheorem{remark}{Remark}[section]

\newtheorem{pro}{Proposition}[section]

\def\bc{\begin{center}}
\def\ec{\end{center}}

\def\hang{\hangindent\parindent}
\def\textindent#1{\indent\llap{\qquad #1\ \ \enspace}\ignorespaces}
\def\ref{\par\hang\textindent}

\begin{document}
 \abovedisplayskip=8pt plus 1pt minus 1pt
\belowdisplayskip=8pt plus 1pt minus 1pt
%-------------------  First Head  -----------------------------------------
\thispagestyle{empty} \vspace*{-3.0truecm} \noindent
\parbox[t]{6truecm}{\footnotesize\baselineskip=11pt\noindent  {} %Acta Mathematica
%%Sinica, English Series\\
%%1999, Jan., Vol.15, No.1, p. 1--11\\
%%Http://www.ActaMath.com\\
%DOI:
 } \hfill
%%\parbox[c]{6truecm}{\vbox{\hsize 3.6576 true cm %
%%  \vskip 3.8 true cm %1.8373
%%  \relax\hbox to0.4\hsize{\hbox to0pt{\special{BMF=actmark.BMF}}\hss}\hss}}
%\hbox to\textwidth{\vbox{\footnotesize\baselineskip=11pt\noindent
%Acta Mathematica Sinica, English Series\hfill LOGO\\
%1999, Jan., Vol.15, No.1, pp. 1--11\hfill \copyright Spring-Verlag 1999}
%\vbox{\hsize3.6576 true cm
%  \vskip1.8373 true cm
%  \relax\hbox to\hsize{\hbox to0pt{\special{BMF=ACTMARK.BMF}}\hss}\hss}}
%===================Text=============================================
\vspace{1 true cm}

\bc{\large\bf Dynkin game under ambiguity in continuous time }
%\footnotetext{\footnotesize Received February 24, 1998, Revised
%September 1, 1998, Accepted September 9, 1996}}
\footnotetext{\footnotesize + Corresponding author.} \ec
\vspace*{0.1 true cm}
\bc{\bf Helin Wu$^{+}$  \\
{\small\it School of Mathematics, Shandong University, Jinan 250100, China\\
\small\it \quad E-mail: wuhewlin@gmail.com}}\ec \vspace*{3 true mm}

\begin{abstract}
In this paper, we want to investigate some kind of Dynkin's game
under ambiguity which is represented by Backward Stochastic
Differential Equation (shortly BSDE) with standard generator
function $g(t,y,z)$. Under regular assumptions, a pair of saddle
point can be obtained and the existence of the value function
follows. The constrained case is also treated in this paper.

\end{abstract}
{\bf Keywords:}Ambiguity, BSDE, Dynkin game, RBSDE.

\section {Introduction}
Dynkin's stopping  games   was first introduced and studied by
Dynkin in [3], and was generalized in  J.Neveu  [9],
N.V.EIbakidze[11],  Yu.I.Kifer [18] , Y.Ohtsubo [19], [20], [21]
etc. with discrete parameter with or without  a finite constraint.
The continuous time version was  also studied in many literature
(for examples,  H.Morimoto [5],  L.Stettner [10] and N.V.Krylov [12]
etc.). We want to investigate some kind of Dynkin's game under
ambiguity in continuous time  in this paper.

A general  formulation of Dynkin's game states  as follows. Define
the lower and upper value function as

$$\underline{V}_t:=ess\sup_{\tau\in{\mathcal{T}_t}}ess\inf_{\sigma\in{\mathcal{T}_t}}E[R_t(\tau,\sigma)|\mathcal{F}_t],\eqno(1.1)$$
and
$$\overline{V}_t:=ess\inf_{\sigma\in{\mathcal{T}_t}}ess\sup_{\tau\in{\mathcal{T}_t}}E[R_t(\tau,\sigma)|\mathcal{F}_t].\eqno(1.2)$$
where $R_t(\tau,\sigma)$ is a function of two stopping times $\tau$
and $\sigma$ satisfying some suitable assumptions. One often try to
find sufficient conditions when $\overline{V}_t =\underline{V}_t$
holds. It is easy to see that $\overline{V}_t \geq\underline{V}_t$,
to get the reverse inequality, one often look for  a pair of saddle
point $(\tau_t^*,\sigma_t^*)$ for which
$$E[R_t(\tau,\sigma_t^*)|\mathcal{F}_t]\leq
E[R_t(\tau_t^*,\sigma_t^*)|\mathcal{F}_t]\leq
E[R_t(\tau_t^*,\sigma)|\mathcal{F}_t]\eqno(1.3)$$ holds for any
$\sigma$ and $\tau$ taking  values in $t$ and $T$. If $(1.3)$ is
true, then $V(t):=\overline{V}_t =\underline{V}_t$ by the definition
of $(1.1)$ and $(1.2)$ and $V(t)$ is called as the value function of
the Dynkin's stopping game.

There are many ways to solve this game. Since stopping game is an
extension of optimal stopping problem, the martingale approach is a
nice choice. In fact, we can find  a pair of saddle point and the
value function   by solving a double optimal stopping problem, for
reference see E.B.Dynkin [3], N.V.Krylov [12] etc. Since Reflected
Backward Stochastic Differential Equation (shortly for RBSDE)  with
lower barrier has been proved  useful  to solve optimal stopping
problem, some author find out a way to solve  Dynkin's game by
solving RBSDE with lower and upper obstacles in J.Cvitanic;
I.Karatzas [8], S.Hamad\`{e}ne; J.-P.Lepeltier [17] etc. Moreover,
A.Bensoussan; A.Friedman [2] and  A.Friedman [1] developed the
analytical theory of stochastic differential games with stopping
times in Markov setting, they studied the value and saddle-points of
such a game using appropriate partial differential equations ,
variational inequalities, and free-boundary problems. Of course,
there are still other ways  to solve this game   such as by pathwise
approach (see I.KaratzasSource [7]) and by connection with singular
control problem ( see  I.Karatzas; H.Wang [6]  ).

Inspired by J.Cvitanic; I.Karatzas [8], in this paper we want to
study a similar Dynkin's game in the stochastic  environment  with
ambiguity and we evaluate the reward process by nonlinear
$g$-expectations. More explicitly, our problem can be formulated as
follows. We define the lower value function and the upper value
function as

$$\underline{V}_t:=ess\sup_{\tau\in{\mathcal{T}_t}}ess\inf_{\sigma\in{\mathcal{T}_t}}\mathcal{E}^g_t[R(\tau,\sigma)]\eqno(1.4)$$
and
$$\overline{V}_t=ess\inf_{\sigma\in{\mathcal{T}_t}}ess\sup_{\tau\in{\mathcal{T}_t}}\mathcal{E}^g_t[R(\tau,\sigma)]\eqno(1.5)$$
respectively.

 Where
$R(\tau,\sigma):=L(\tau)1_{(\tau\leq\sigma)}+U(\sigma)1_{(\sigma<\tau)}$
and $\mathcal{T}_t$ are stopping times taking values  between $t$
and $T$, the finite termination of problem. Under some suitable
assumptions on the two processes $L(t)$ and $U(t)$,  we want to find
out a pair of saddle point $(\tau_t^*,\sigma_t^*)$ such that
$$\mathcal{E}^g_t[R(\tau,\sigma_t^*)]\leq
\mathcal{E}^g_t[R(\tau_t^*,\sigma_t^*)]\leq
\mathcal{E}^g_t[R(\tau_t^*,\sigma)]\eqno(1.6)$$ for any
$\tau,\sigma\in \mathcal{T}_t$ and  then by definition of $(1.4)$
and $(1.5)$, the game has a value function.

This problem looks very like with  the problem  stated and solved in
J.Cvitanic; I.Karatzas [8] , but there is difference between them,
although the solutions are  same as we will prove. To make our
problem meaningful, we will point out the main difference in section
3 and treat a more complicated case with constraints, in which case
we evaluate reward process by $g_\Gamma$-expectation whose
definition will be given in section 2.

 Our paper organized as follows.
 In section two,  the  necessary  framework  and some useful propositions of BSDE was reviewed,
 and the main result and its proof is stated in section 3.
\section {BSDE, Reflected BSDE and Constrained BSDE}
Given a probability space $(\Omega,\mathcal{F},P)$ and $R^d$-valued
Brownian motion \mbox{$W(t)$}, we consider a sequence
$\{(\mathcal{F}_t);t\in[0,T]\}$ of filtrations generated by Brownian
motion $W(t)$  and augmented by P-null sets. $\mathcal{P}$  is the
$\sigma$-field of predictable sets of  $\Omega\times[0,T]$. We use
$L^2(\mathcal{F}_T)$ to denote the space of all $F_T$-measurable
random variables $\xi:\Omega\rightarrow R^d$ for which
$$\parallel \xi\parallel^2=E[|\xi|^2]<+\infty.$$
and use $H_T^2(R^d)$ to denote  the space of predictable process
$\varphi:\Omega\times[0,T]\rightarrow R^d$ for which
$$\parallel \varphi\parallel^2=E[\int_0^T|\varphi|^2]<+\infty.$$

$S^k_n$ denotes the space of $(\mathcal{F}_t)$-progressively
measurable processes $\varphi: [0,T]\times\Omega\mapsto R^n$ with
$E(\sup_{0\leq t\leq T}\parallel \varphi\parallel^k)<\infty, k\in
N.$

 $S^2_{ci} $ denotes the space of continuous, increasing, $(\mathcal{F}_t)$-adapted processes $A:
 [0,T]\times \Omega\mapsto [0, \infty)
$ with $A(0)=0, E(A^2(T)) < \infty$.

Given a function $g:[0,T]\times R\times R^d\rightarrow R$, following
assumptions always used in theory of BSDE.
$$|g(\omega,t,x_1,y_1)-g(\omega,t,x_2,y_2)|\leq M(|x_1-x_2|+|y_1-y_2|), \ \ \forall
(x_1,y_1),(x_2,y_2)\eqno(A1)$$ for some $M>0$.

$$g(\cdot,x,y)\in H_T^2(R)\quad \forall x\in R,\, y\in R^d.  \eqno(A2)$$

The  BSDE  driven by $g(t,x,y)$ is given by
$$-dX(t)=g(t,X(t),Y(t))dt-Y'(t)dW(t) \eqno(2.1)$$
 where $X(t)\in R$ and $W(t)\in R^d$.
Suppose that $\xi\in L^2(\mathcal{F}_T)$ and $g$ satisfies (A1) and
(A2), E.Pardoux  and  S.G.Peng[4] proved the existence  of adapted
solution $(X(t),Y(t))$ of such BSDE. We call $(g,\xi)$ standard
parameters for the BSDE.

For later use, we collect some useful propositions of BSDE below,
its proof can be found in many papers such as  S.G.Peng[13].
\begin{pro}
If the generator function $g(t,x,y):[0,T]\times R\times R^d\mapsto
R$ satisfies assumptions (A1) and (A2).
 For any stopping time
$\tau\leq\sigma\leq T$,  we denote $X(\tau)$ in the following BSDE
as $\mathcal{E}^g_{\tau,\sigma}(\zeta)$ ,
$$X(\tau)=\zeta+\int_\tau^\sigma g(s,X(s),Y(s))ds -\int_\tau^\sigma Y'(s)dW(s).$$
where $\zeta$ is $\mathcal{F}_\sigma$ measurable, then we have
\begin{itemize}
\item[(i)](Comparison proposition) If $\mathcal{F}(\sigma)$-measurable variables $\xi \geq \eta \,a.s$ , then
$\mathcal{E}^g_{\tau,\sigma}(\xi)\geq
\mathcal{E}^g_{\tau,\sigma}(\eta)$ for any stopping times
$0\leq\tau\leq\sigma\leq T$ a.s.
\item[(ii)] If $g(t,y,0)=0$, then for any stopping times
$\tau\leq\sigma\leq T$ and  $\mathcal{F}_\sigma$-measurable variable
$\zeta$, we have
$\mathcal{E}^g_{\tau,\sigma}(\zeta)=\mathcal{E}^g_{\tau,T}(\zeta)$
and we write $\mathcal{E}^g_{\tau,\sigma}(\zeta)$ shortly as
$\mathcal{E}^g_\tau(\zeta)$ when $\sigma=T$.
\item[(iii)] (Coherence) If $g(t,y,0)=0$, then for any stopping
times $\tau\leq\sigma\leq T$ and $\mathcal{F}_\sigma$-measurable
variable $\zeta$, we have
$\mathcal{E}^g_\tau(\mathcal{E}^g_\sigma(\zeta))=\mathcal{E}^g_\tau(\zeta).$
\end{itemize}
\end{pro}

The theory of BSDE has wildly used in many fields such as financial
mathematics and stochastic optimal control problems.  Some brilliant
use of these is that one can connect the optimal stopping problem
with BSDE reflected by some lower barrier and connect Dynkin's game
problem with BSDE reflected from below and above  by lower barrier
and upper barrier respectively. Here, for later  proof  of our
problem we need to introduce Reflected  BSDE.

\begin{definition}
 (Backward stochastic differential equation (BSDE) with
upper and lower reflecting barriers). Let $\xi$ be a given random
variable in $L^2(\mathcal{F}_T)$, and $g:[0,T]\times\Omega\times R
\times R^d \mapsto R$ a given $\mathcal{P}\otimes
\mathcal{B}(R)\otimes \mathcal{B}(R^d)$-measurable functions
satisfying $(A1)$ and $(A2)$.

Consider  two continuous processes $L,U$ in $S^2_1$  satisfying
$$L(t)\leq U(t), \qquad \forall\, 0\leq t \leq T \quad \text{and}
\quad L(T)\leq\xi\leq U(T) \quad a.s.$$

 We say that a triple $(X,Y,K)$ of F-progressively measurable
processes  $X: [0,T]\times \Omega \mapsto R$, $Y: [0,T]\times \Omega
\mapsto R^d$  and $K: [0,T]\times \Omega \mapsto R$ is a solution of
the Backward Stochastic Differential Equation (BSDE) with reflecting
barriers $U, L$ (upper and lower, respectively), terminal condition
$\xi$ and coefficient $g$, if the following hold:

\begin{itemize}
\item[(i)] $K=K^+-K^-$ with $K^{\pm}\in S^2_{ci}$,

\item[(ii)] $Y\in H^2_d$,
\end{itemize}
and

$$X(t) = \xi+\int_t^Tg(s,X(s),Y(s))ds
+K^+(T)-K^+(t) \\
-(K^-(T)-K^-(t))-\int_t^TY'(s)dW(s), \qquad 0\leq t \leq
T,\eqno(2.2)$$

$$L(t) \leq X(t) \leq U(t), \qquad 0\leq t \leq T,\eqno(2.3)$$

$$\int_0^T(X(t)-L(t))dK^+(t)=\int_0^T(U(t)-X(t))dK^-(t)=0,\eqno(2.4)$$

almost surely.

\end{definition}

The processes $L, U$ play the role of reflecting barriers, these are
allowed to be random and time-varying, and the state-process $X$ is
not allowed to cross them  on its way to the prescribed terminal
target condition $X_T=\xi$. J.Cvitanic; I.Karatzas [8] has solved
this kind of BSDE with two reflected barriers in $S_1^2$  via
solving Dynkin'game and double optimal stopping problem.
 S.G.Peng and M.Y.Xu [16] treated different case  of such problem with different class of
barriers which is sufficient for the use  of our paper.

We will also treat constrained case  in our paper, so we introduce
Constrained Backward Stochastic Differential Equation (CBSDE) at the
same time.

The constraints in our paper is like that in S.G.Peng [14], namely
for a given nonnegative function $\phi(t,x,y): [0,T]\times R\times
R^d\mapsto R^+$  we ask the solution $(X(t),Y(t))$ of BSDE
satisfying
$$\phi(t,X(t),Y(t))=0, a.s \quad\text{for any} \quad t\in[0,T].\eqno(C)$$

In constrained case, it often need an increasing process acting as
singular control to force  the solution stays in the constrained
filed. BSDE with an increasing process is called a
$g$-super-solution and the smallest one plays a crucial role.
\begin{definition}
(super-solution) A super-solution of a BSDE associated with the
standard parameters $(g,\xi)$ is a vector process $(X(t),Y(t),C(t))$
satisfying
$$-dX(t)=g(t,X(t),Y(t))dt+dC(t)-Y'(t)dW(t),\quad X(T)=\xi,\eqno(2.5)$$
or being equivalent to
$$X(t)=\xi+\int_t^Tg(s,X(s),Y(s))ds-\int_t^TY'(s)dW(s)+\int_t^TdC(s), \eqno(2.5')$$
where $(C_t,t\in[0,T])$ is an increasing, adapted, right-continuous
process with $C_0=0$.

\end{definition}

\begin{definition}( $g_\Gamma$-solution or the minimal solution ) A g-super-solution $(X(t), Y(t), C(t))$ is said to be the
 the minimal solution of a constrained backward differential stochastic equation (shortly CBSDE), given $y_T=\xi$,
 subjected to the constraint $(C)$ if for any other g-super-solution
$(\tilde{X}(t),\tilde{Y}(t), \tilde{C}(t))$ satisfying $(C)$  with
$\tilde{X}(T)=\xi$, we have $X(t)\leq\tilde{X}(t) $ a.e., a.s.. When
both $g(t,x,0)=0$ and $\phi(t,x,0)=0$ for any $x\in R,t\in [0,T]$,
the minimal solution is denoted by $\mathcal{E}_t^{g,\phi}(\xi)$ and
for convenience called as $g_\Gamma$-solution. Sometimes, we also
call $g_\Gamma$-expectation $\mathcal{E}_t^{g,\phi}(\xi)\triangleq
X(t)$ the dynamic $g_\Gamma$-expectation with constraints $(C)$.
\end{definition}
For any $\xi\in L^2(\mathcal{F}_T)$, we denote
$\mathcal{H}^{\phi}(\xi)$ as the set of g-super-solutions
$(X(t),Y(t),C(t))$ subjecting to $(C)$ with $X(T)=\xi$.  When
$\mathcal{H}^{\phi}(\xi)$ is not empty,   S.G.Peng [14] proved that
$g_\Gamma$-solution exists.

Similarly with S.G.Peng; M.Y.Xu [15], let
$$\Gamma_t:=\{(t,x,y)\in [0,T]\times R\times
R^d|\phi(t,x,y)=0\},$$

the  $g_\Gamma$-solution is defined as the smallest
$g$-super-solution  with constraints (C).

We give a continuous property of $g_\Gamma$-solution for later use.
\begin{pro}
Suppose  the generator function $g(t,x,y)$ and the  constraint
function $\phi(t,x,y)$ both satisfy conditions (A1) and (A2),
  $\{\xi_n\in L_T^2(P), \ n=1,2,\cdots\}$ is an bounded increasing
sequence and converges almost surely to $\xi\in L_T^2(P)$, if
$\mathcal{E}_t^{g,\phi}(\zeta)$ exists for
$\zeta=\xi,\xi_n,n=1,2,\cdots$, then
$$\lim_{n\rightarrow \infty}\mathcal{E}_t^{g,\phi}(\xi_n)= \mathcal{E}_t^{g,\phi}(\xi) \quad a.s \quad\forall \,t\in [0,T]. $$
\end{pro}

\proof Without noting, all the proofs  go on   under almost surely.

According to S.G.Peng [14], the solutions  $x^m(t)(\xi)$ of
$$x^m(t)(\xi)=\xi + \int_t^Tg(x^m(s)(\xi),y^m(s),s)ds+A^m(T)-A^m(t)-\int_t^Ty^m(s)dW(s).$$
is an increasing sequence and converges  to
$\mathcal{E}_t^{g,\phi}(\xi)$, where $$A^m(t): =
m\int_0^t\phi(x^m(s),y^m(s),s)ds.$$

It is easy to see  $\{\mathcal{E}_t^{g,\phi}(\xi_n),n=1,2,\cdots\}$
is an increasing sequence. We denote its limit at $t$ as  $a_t$,
then $a_t\leq \mathcal{E}_t^{g,\phi}(\xi)$. Since $\xi_n$ converges
almost surely increasingly to $\xi\in L_T^2(R)$, by dominated
convergence theorem, it also converges strongly  in $L_T^2(R)$, then
by the continuous dependence property of g-supersolution, the limit
of $\{x^m(t)(\xi_n)\}_{n=1}^{\infty}$ is $x^m(t)(\xi)$ for any fixed
$m$.

We want to show  that $a_t=\mathcal{E}_t^{g,\phi}(\xi)$. If on the
contrary on has  $a_t<\mathcal{E}_t^{g,\phi}(\xi)$, then there is
some $\delta>0$ such that
$\mathcal{E}_t^{g,\phi}(\xi)-\mathcal{E}_t^{g,\phi}(\xi_n)>\delta$
for any $n$. On the other hand, for any $\epsilon>0$, $0\leq
\mathcal{E}_t^{g,\phi}(\xi)-x^m(t)(\xi)\leq\epsilon$ holds for some
larger $m_0$. Fixing  $m_0$, $\epsilon$, there is some $n_0$ which
depends on $m_0$ and $\epsilon$ such that $0\leq x^{m_0}(t)(\xi)-
x^{m_0}(t)(\xi_{n_0})\leq\epsilon$, so
$\mathcal{E}_t^{g,\phi}(\xi)-x^{m_0}(t)(\xi_{n_0})\leq 2\epsilon$,
but we have $\mathcal{E}_t^{g,\phi}(\xi)-x^{m_0}(t)(\xi_{n_0})\geq
\mathcal{E}_t^{g,\phi}(\xi)-\mathcal{E}_t^{g,\phi}(\xi_{n_0})>\delta$,
this is impossible for $\epsilon<\frac{\delta}{2}$.
\hspace*{\fill}$\Box$

\section {Dynkin's game under ambiguity }

In this section we firs review some  exited result  about Reflected
BSDE and Dykin' game.

In  [8], if  $(X,Y,Z)$ is the solution of Reflected BSDE stated in
above section, then it is said that $X(t)$ equals the value function
of the   Dynkin's game of (1.1) and (1.2) with
$$R_t(\tau,\sigma)=\int_t^{\tau\wedge \sigma}g(s,(X(s),Y(s))ds +
L(\tau)1_{(\tau<T,\tau\leq\sigma)}+U(\sigma)1_{(\sigma<\tau)}+\xi
1_{\tau\wedge\sigma=T}.\eqno(I)$$

More generally in S.Hamad\`{e}ne, J.-P.Lepeltier [17], the author
considered the mixed zero-sum stochastic differential game with
payoff
$$J(u,\tau;v,\sigma)=E^{(u,v)}[\int_0^{\tau\wedge\sigma}g(s,X(s),u(s),v(s))ds
+ L(\tau)1_{\tau\leq\sigma,\sigma<T}+ U(\sigma)1_{\sigma<\tau}+\xi
1_{\tau\wedge\sigma=T}].\eqno(II)$$

Under the assumption  $g(t,x,0)=0$,
 we can  explore
$\mathcal{E}^g_t[L(\tau)1_{(\tau\leq\sigma)}+U(\sigma)1_{(\sigma<\tau)}]$
as
$$E[\int_t^{\tau\wedge \sigma}g(s,(X^{\tau,\sigma}(s),Y^{\tau,\sigma}(s))ds +
L(\tau)1_{(\tau<T,\tau\leq\sigma)}+U(\sigma)1_{(\sigma<\tau)}+\xi
1_{\tau\wedge\sigma=T}|\mathcal{F}_t]\eqno(III)$$ with $\xi=L(T)$.

From the expression of (I), (II) and (III) above, we can easily find
that the difference between integrands  in the integral. In (I),
$(X(s),Y(s))$ is fixed at first, in (II), $X(s)$ only depends on
controls $(u,v)$, but in our problem  (III)
$((X^{\tau,\sigma}(s),Y^{\tau,\sigma})$ depends on stopping times
$(\tau,\sigma)$.

But with the help of Reflected BSDE and comparison proposition  of
BSDEs or $g$-martingale theory, we can find a pair of saddle point
of such Dynkin' game under nonlinear expectation. Below is our main
result in unconstrained case.
\begin{thm}
Let $L(t),U(t)\in S_1^2$ and $L(t)\leq U(t), 0\leq t\leq T$. The
generator function $g(t,x,y)$ satisfies assumptions $(A1),A(2)$ and
$g(t,x,0)=0,\forall t,x$. Suppose $(X(t), Y(t),K(t))$ is the
solution of Reflected BSDE  as formulated in definition $(2.1)$ with
terminal value $L(T)$, then the Dynkin's game stated in
$(1.4),(1.5)$ has a pair of saddle point $(\tau_t^*,\sigma_t^*)$ and
hence the value function exists. Furthermore, the pair of saddle
point can  be represented by
$$\tau_t^*=\inf\{s\geq t: L(s)=X(s)\}\wedge T\eqno(3.1)$$
and
$$\sigma_t^*=\inf\{s\geq t:U(s)=X(s)\}\wedge T.\eqno(3.2)$$
and
$$\underline{V}(t)=\overline{V}(t)=X(t).$$
\end{thm}
\proof We want to prove $(1.6)$.

Fix $\sigma_t^*$ first and let $\tau$ be arbitrary stopping time
taking values in $[t,T]$, then we have
\begin{eqnarray*}
\mathcal{E}^g_t[R_t(\tau,\sigma_t^*)]&=&\mathcal{E}^g_t[L(\tau)1_{(\tau\leq\sigma_t^*)}+U(\sigma_t^*)1_{(\sigma_t^*<\tau)}]\\
&\leq&\mathcal{E}^g_t[X(\tau)1_{(\tau\leq\sigma_t^*)}+U(\sigma_t^*)1_{(\sigma_t^*<\tau)}]
  \quad(\text{for} \,L\leq X \,\text{and}\, X(\sigma_t^*)= U(\sigma_t^*))\\
  &=&\mathcal{E}^g_t[X(\tau\wedge \sigma_t^*)].
\end{eqnarray*}

At this step, we need to prove that
$$\mathcal{E}^g_t[X(\tau\wedge \sigma_t^*)]\leq \mathcal{E}^g_t[X(\tau_t^*\wedge
\sigma_t^*)]\eqno(3.3)$$
 for any $\tau$ values in $[t,T]$.
 By $(2.4)$ and $(3.1),(3.2)$, we have that when $\tau\leq\tau_t^*$,
 $A^+(\tau)=A^+(t)$; when $\sigma\leq\sigma_t^*$,
 $A^-(\sigma)=A^-(t)$.

 So on the set $(\tau\leq\tau_t^*)$, by the equation of (2.2), we
 have
 $$X(\tau\wedge\sigma_t^*) = X(\tau_t^*\wedge\sigma_t^*)+\int_{\tau\wedge\sigma_t^*}^{\tau_t^*\wedge\sigma_t^*}g(s,X(s),Y(s))ds
-\int_{\tau\wedge\sigma_t^*}^TY'(s)dW(s),$$

this means

$$X(\tau\wedge \sigma_t^*)=
\mathcal{E}^g_{\tau\wedge \sigma_t^*}[X(\tau_t^*\wedge \sigma_t^*)].
\eqno(3.4)$$

On the other hand, when $\tau>\tau_t^*$, similarly we have
$$X(\tau_t^*\wedge\sigma_t^*) =
X(\tau\wedge\sigma_t^*)+\int_{\tau_t^*\wedge\sigma_t^*}^{\tau\wedge\sigma_t^*}g(s,X(s),Y(s))ds
+A^+(\tau\wedge\sigma_t^*)-A^+(\tau\wedge\sigma_t^*)
-\int_{\tau_t^*\wedge\sigma_t^*}^{\tau\wedge\sigma_t^*}Y'(s)dW(s),$$
and this means
$$\mathcal{E}^g_{\tau_t^*\wedge \sigma_t^*}[X(\tau\wedge
\sigma_t^*)]\leq X(\tau_t^*\wedge \sigma_t^*). \eqno(3.5)$$

Taking  conditional $g$-expectation  on both hands of (3.4) and
(3.5), by the coherence property of $g$-solutions, (iii) of
proposition 2.1,  one has (3.3).

Now, we fix $\tau_t^*$ , then for any $\sigma$ taking values in
$[t,T]$, we want to show that

$$\mathcal{E}^g_t[X(\tau_t^*\wedge \sigma_t^*)]\leq \mathcal{E}^g_t[X(\tau_t^*\wedge
\sigma)].\eqno(3.6)$$

Similarly, when $(\sigma\leq\sigma_t^*)$,  we have

$$X(\tau_t^*\wedge\sigma) = X(\tau_t^*\wedge\sigma_t^*)+\int_{\tau_t^*\wedge\sigma}^{\tau_t^*\wedge\sigma_t^*}g(s,X(s),Y(s))ds
-\int_{\tau_t^*\wedge\sigma}^{\tau_t^*\wedge\sigma_t^*}Y'(s)dW(s),$$
and
$$X(\tau_t^*\wedge \sigma)=
\mathcal{E}^g_{\tau_t^*\wedge \sigma}[X(\tau_t^*\wedge \sigma_t^*)].
\eqno(3.7)$$ becomes true.

When $\sigma>\sigma_t^*$, (2.2) becomes
$$X(\tau_t^*\wedge\sigma_t^*) =
X(\tau_t^*\wedge\sigma)+\int_{\tau_t^*\wedge\sigma_t^*}^{\tau_t^*\wedge\sigma}g(s,X(s),Y(s))ds
-(A^-(\tau_t^*\wedge\sigma)-A^-(\tau_t^*\wedge\sigma_t^*))
-\int_{\tau_t^*\wedge\sigma_t^*}^{\tau_t^*\wedge\sigma}Y'(s)dW(s)$$
and this means
$$\mathcal{E}^g_{\tau_t^*\wedge \sigma_t^*}[X(\tau_t^*\wedge
\sigma)]\geq X(\tau_t^*\wedge \sigma_t^*). \eqno(3.8)$$

Taking $g$-expectation in (3.7) and (3.8), we
have(3.6).\hspace*{\fill}$\Box$

\begin{remark}
 Comparing  the method used   to prove  Theorem 4.1   in J.Cvitanic; I.Karatzas [8] and the
 method to prove 3.1 in our paper, although they are very similar, but the
 advantages of our problem helps us to use $g$-martingale theory
 directly and this is very convenient for us to handle constrained
 case later.
\end{remark}

 Next, we  then go to the constrained case, that is we evaluate the reward
 process by the constrained $g$-expectation which is also named as
 $g_\Gamma$-expectation in  S.G.Peng; M.Y.Xu [15].
 Define similarly lower and upper value function as non-constrained
 case,
$$\underline{V}_t:=ess\sup_{\tau\in{\mathcal{T}_t}}ess\inf_{\sigma\in{\mathcal{T}_t}}\mathcal{E}^{g,\phi}_t[R(\tau,\sigma)]\eqno(3.9)$$

$$\overline{V}_t=ess\inf_{\sigma\in{\mathcal{T}_t}}ess\sup_{\tau\in{\mathcal{T}_t}}\mathcal{E}^{g,\phi}_t[R(\tau,\sigma)]\eqno(3.10)$$
with
$R(\tau,\sigma)=L(\tau)1_{(\tau\leq\sigma)}+U(\sigma)1_{(\sigma<\tau)}.$

 For any integer $m$,  let $g_m=g+m\phi$. For any  fixed $m$,
there is an unique solution of  Reflected $g_m$-solution with double
barriers (from below and above by $L$ and $U$ respectively) and we
denote it as $(X^m,Y^m,K^m), m=1,2,\cdots$. Since $X^m$ can be
obtained by penalization  method, the comparison proposition ensures
that  $\{X^m\}$ is an increasing sequence of process.

Define the corresponding pairs of  stopping  times as
$$\tau_t^*(m)=\inf\{s\geq t: L(s)=X^m(s)\}\wedge T\eqno(3.11)$$
and
$$\sigma_t^*(m)=\inf\{s\geq t:U(s)=X^m(s)\}\wedge T.\eqno(3.12)$$

It is easy to see that $\{\tau_t^*(m)\}$  is increasing and
$\{\sigma_t^*(m)\}$ is decreasing and
$$\tau_t^*=\lim_{m\rightarrow \infty}\tau_t^*(m),\,\sigma_t^*=\lim_{m\rightarrow
\infty}\sigma_t^*(m)\eqno(3.13)$$ are stopping times.

With these in hand, we can then state and prove  our result in
constrained case below.
\begin{thm}
Let $g$ and $\phi$ satisfy assumptions (A1) and (A2), $L(t)$ and
$U(t)$ are    nonnegative  continuous processes and there is some
constant $B>0$ such that $L(t)\leq B, U(t)\leq B $ for any $t\in
[0,T]$. We consider the Dynkin's game with lower and upper value
function defined in (3.9) and (3.10). If $L(t)$ is increasing, then
the pair of stopping time defined in (3.13) is a saddle point.
\end{thm}
\proof For any  $n\leq m$, by comparison theorem of BSDE and results
obtained in unconstrained case of Theorem 3.1, we have
$$\mathcal{E}^{g_n}_t[X(\tau\wedge \sigma_t^*(m))]\leq\mathcal{E}^{g_m}_t[X(\tau\wedge \sigma_t^*(m))]
\leq\mathcal{E}^{g_m}_t[X(\tau_t^*(m)\wedge
\sigma_t^*(m))]=X^m(t)\eqno(3.14)$$ for any $\tau$ and
$\tau_t^*(m),\sigma_t^*(m)$ defined  in (3.11) and (3.12).

On the other hand, one has
$$X^m(t)=\mathcal{E}^{g_m}_t[X(\tau_t^*(m)\wedge \sigma_t^*(m))]\leq\mathcal{E}^{g_m}_t[X(\tau_t^*(m)\wedge \sigma)]
\leq\mathcal{E}^{g,\phi}_t[X(\tau_t^*(m)\wedge \sigma)]\eqno(3.15)$$
for any $\sigma$ taking values in $[0,T]$.

 Firs, we  take limit in  (3.14) and (3.15) as $m\rightarrow
 \infty$,  set $X(t):=\lim_{m\rightarrow\infty}X^m(t)$,
we have
$$\mathcal{E}^{g_n}_t[X(\tau\wedge \sigma_t^*]\leq X(t)\eqno
(3.16)$$ and
$$X(t)\leq \mathcal{E}^{g,\phi}_t[X(\tau_t^*\wedge
\sigma)].\eqno(3.17)$$

 Since  $g$-solution is continuous dependence on its
terminal value and, by Proposition 2.2, $g_\Gamma$-solution is
continuous from below with its terminal variable, see  also Remark
3.2 in this paper.

By (3.16) and (3.17), we conclude that The Dynkin's game has a value
function X(t).

To prove $(\tau_t^*,\sigma_t^*)$ is a saddle point, we want to prove
$$\mathcal{E}^{g,\phi}_t[X(\tau_t^*\wedge \sigma_t^*)]=X(t)\eqno(3.18)$$.

First note that
$$X^m(t)=\mathcal{E}^{g_m}_t[X(\tau_t^*(m)\wedge \sigma_t^*(m))]\leq\mathcal{E}^{g_m}_t[X(\tau_t^*(m)\wedge \sigma_t^*)]
\leq\mathcal{E}^{g,\phi}_t[X(\tau_t^*(m)\wedge \sigma_t^*)].$$

 Taking limit in above equation  as $m\rightarrow \infty$, because $L(t)$
 is continuous  and increasing and
 $g_\Gamma$-solution is continuous from below, we have
 $$X(t)\leq\mathcal{E}^{g,\phi}_t[X(\tau_t^*\wedge \sigma_t^*)].\eqno(3.19)$$

 For the other side inequality,
note that $\sigma_t^*\leq\sigma_t^*(m)$ for any $m$, we have

$$X^m(t)=\mathcal{E}^{g_m}_t[X(\tau_t^*(m)\wedge \sigma_t^*(m))]=\mathcal{E}^{g_m}_t[X(\tau_t^*(m)\wedge
\sigma_t^*)]\geq\mathcal{E}^{g_m}_t[X(\tau_t^*\wedge \sigma_t^*)].$$
Taking  limit   as $m\rightarrow \infty$ in above equation, on has

$$X(t)\geq\mathcal{E}^{g,\phi}_t[X(\tau_t^*\wedge
\sigma_t^*)].\eqno(3.20)$$

Combine  (3.19) and (3.20) together, we obtain (3.18) and thus
complete our proof.\hspace*{\fill}$\Box$

\begin{remark}
Note that under the assumptions of Theorem 3.2,
$\mathcal{E}^{g,\phi}_t[X(\tau)]$ is meaningful for any stopping
time $\tau$ taking values in $[t,T]$.  It is easy to see that
$X^m(t)\leq B$ for any $t\in [0,T]$ and assumptions (A1) and (A2)
together with $g(t,x,0)=0$ and $\phi(t,x,0)=0$ ensures that
$g_\Gamma$-solution  is well defined on $L_T^\infty(P)$, the space
of  all essentially bounded $\mathcal{F}_T$-measurable  variables.
In the paper S.G.Peng and M.Y.Xu [15], the author defined a new
subspace of $L^2_T(P)$:
$$L^2_{+,\infty}(\mathcal{F}_T)\triangleq\{\xi\in
L^2(\mathcal{F}_T), \xi^+\in L^\infty(\mathcal{F}_T)\}.$$

For any $\xi\in L^2_{+,\infty}(\mathcal{F}_T)$ with terminal
condition $y_T=\xi$,   $g_\Gamma$-solution  exists if
$$g(t,y,0)\leq L_0 +M|y| \quad \text{and}\quad (y,0)\in \Gamma_t\eqno(3.21)$$
holds for a large constant $L_0$ and for any $ y\geq L_0$ and if
there exists a deterministic process $a(t)$ such that $L(t) \leq
a(t)$ on $[0,T ]$. Under assumptions on $g$ and $\phi$ as above
mentioned, (3.21) is satisfied  for any $L_0\geq 0$ and $M$ in
$(A1)$ and we can chose  $a(t)=B$. It is obvious
$L^\infty(\mathcal{F}_T)\subset L^2_{+,\infty}(\mathcal{F}_T)$ and
$g_\Gamma$-solution is defined well on the whole space
$L^\infty(\mathcal{F}_T)$

\end{remark}

\begin{remark}
Roughly speaking, the continuous property from below of
$g_\Gamma$-expectation is a simple consequence  of the fact that we
can change the order of limits in
$\lim_{n\rightarrow\infty}\lim_{m\rightarrow\infty}a_{m,n}$ and
$\lim_{m\rightarrow\infty}\lim_{n\rightarrow\infty}a_{m,n}$  when
$a_{m,n}$ are both increasing with $n$ and $m$. The complete
continuous property is more complicated since it concerns mini-max
problem. But except $g_\Gamma$-expectation is continuous from below,
it is still semi-lower-continuous, and we can still conclude some
useful continuous property by convex assumptions on coefficients of
CBSDE with  help of some wonderful results in convex analysis. Of
course, we can then make different assumptions on $L$ and $U$ to get
through the proof via continuous property in the Theorem 3.2.

\end{remark}

\begin{remark}
It is an open problem  that whether  $X(t)$,  the limit of
$\{X^m(t)\}$ which is a sequence of  Reflected solution of BSDEs, is
still some solution of some kind of Reflected BSDEs?
\end{remark}
\begin{remark}
Dynkin's game problem is very similar with Optimal stopping problem
under ambiguity. It is well known that the Snell envelope  of
Optimal stopping problem for $g$-expectation is same with the
solution of Reflected BSDE with one lower obstacle.  Roughly
speaking, this result  is based on the  following three deep facts:
\begin{itemize}
\item [(i)] The solution of Reflected BSDE with one lower barrier is
the same with the $g_\Gamma$-solution taking the barrier as a
constraint.

\item [(ii)] The Snell envelope  of the barrier under $g$-expectation
is the smallest $g$-super-martingale above the barrier.
\item [(iii)] Under suitable assumptions, $g$-super-solution is  equivalent to
$g$-super-martingale.
\end{itemize}
But in game case, there are no corresponding  theories, so the
results in this paper are not  simple extensions of Optimal stopping
case.
\end{remark}

{}

\end{document}